\documentclass[12pt,reqno]{amsart}

\usepackage{amssymb}
\usepackage{amsmath}
\usepackage{amsthm}
\usepackage{mathrsfs}

\numberwithin{equation}{section}

\allowdisplaybreaks

\theoremstyle{plain}
\newtheorem{thm}{Theorem}[section]
\newtheorem{prop}[thm]{Proposition}
\newtheorem{lem}[thm]{Lemma}

\theoremstyle{definition}
\newtheorem{df}{Definition}[section]

\DeclareMathOperator{\an}{and}

\DeclareMathOperator{\re}{Re}

\DeclareMathOperator{\supp}{supp}

\def\a{\alpha}

\def\bd{D_t}
\def\cd{\d_t}
\def\d{\partial}

\def\dis{\displaystyle}

\def\ov{\overline}

\def\tg{\tilde g}
\def\th{\theta}
\def\tm{\times}

\def\wt{\widetilde}

\def\C{\mathbb{C}}
\def\D{\mathcal{D}}

\def\G{\Gamma}
\def\H{\mathcal{H}}

\def\Lc{\mathcal{L}}

\def\N{\mathbb{N}}
\def\O{\Omega}
\def\Ov{\overline{\O}}
\def\R{\mathbb{R}}

\def\Sg{\Sigma}

\def\smij{\sum_{i,j=1}^d}

\def\({\left(}
\def\){\right)}
\def\<{\left\langle}
\def\>{\right\rangle}

\title[Non-homogeneous boundary value problems]
{Non-homogeneous boundary value problems for fractional diffusion equations in 
$L^2$-setting}

\begin{document}

\begin{abstract}
  In the present article, we study the diffusion equations with fractional 
time derivatives.
  The aim of this paper is to investigate the best possible regularity for the 
initial value/boundary value problems with non-homogeneous Dirichlet boundary 
data.
  The main tool we use here is called the transposition method.
\end{abstract}

\author{Kenichi Fujishiro}



\thanks{Department of Mathematical Sciences, University of Tokyo, 
Komaba, Meguro, Tokyo 153, Japan}

\email{kenichi@ms.u-tokyo.ac.jp}

\date{}

\maketitle

\section{Introduction}\label{sec:int}

  Let $\O$ be a bounded domain of $\R^d$ with $C^2$ boundary $\G:=\d\O$ and 
set $Q:=\O\tm(0,T)$ and $\Sg:=\G\tm(0,T)$.
  We consider the following initial value/boundary value problem for a partial
differential equation with the fractional derivative in time $t$:
\begin{equation}\label{eq:nonhomo}
\begin{cases}
	\cd^{\a}u+A u=0&\mbox{in}\ Q, \\
	u=g		&\mbox{on}\ \Sg, \\
	u(\cdot,0)=0	&\mbox{in}\ \O
\end{cases}
\end{equation}
  with $0<\a<1$.
  Here $\cd^{\a}$ denotes the Caputo derivative, which is defined by
\[
	\cd^{\a}u(x,t)=
	\frac{1}{\G(1-\a)}\int_0^t(t-\tau)^{-\a}
	\frac{\d u}{\d\tau}(x,\tau)d\tau,
\]
  and $\G(\cdot)$ is the Gamma function (see Podlubny \cite{pd}).
  The differential operator $A$ is given by
\[
	A u(x)=
	-\smij\frac{\d}{\d x_i}\(a_{ij}(x)\frac{\d u}{\d x_j}(x)\)+c(x)u(x),
	\quad x\in\O,
\]
  and the coefficients satisfy the following:
\begin{align*}
	&a_{ij}=a_{ji}\in C^1(\Ov),\ 1\le i,j\le d,\quad
		\smij a_{ij}(x)\xi_i\xi_j\ge\mu|\xi|^2,
		\ x\in\Ov,\ \xi\in\R^d, \\
	&c\in C(\Ov),\quad c(x)\ge0,\ x\in\Ov,
\end{align*}
  where $\mu>0$ is constant.
  The function $g$ is given on $\Sg$.

  In the present paper, we study the regularity of the solution to 
\eqref{eq:nonhomo} in the sense of Sobolev spaces.
  As for this problems, Lions and Magenes \cite{lm} showed the result for the 
parabolic equations.

\section{Main Result}\label{sec:main}

  In this section, we prepare the notations and state our main results.
  We denote by $H^s(\O)$, $s\ge0$, the Sobolev spaces.
  For $r,s\ge0$, we abbreviately set
\[
	H^{r,s}(Q):=L^2(0,T;H^r(\O))\cap H^s(0,T;L^2(\O)).
\]
  Then $H^{r,s}(Q)$ is a Hilbert space with the norm $\|\cdot\|_{H^{r,s}(Q)}$ 
given by
\[
	\|u\|_{H^{r,s}(Q)}^2
	:=\|u\|_{L^2(0,T;H^r(\O))}^2+\|u\|_{H^s(0,T;L^2(\O))}^2.
\]
  In particular,
\[
	H^{0,0}(Q)=L^2(Q)=L^2(0,T;L^2(\O)).
\]
  Similarly we set
\begin{align*}
	H^{r,s}(\Sg):=L^2(0,T;H^r(\G))\cap H^s(0,T;L^2(\G)), \\
	H^{r,s}_{\ ,0}(\Sg):=L^2(0,T;H^r(\G))\cap H^s_0(0,T;L^2(\G)).
\end{align*}
  with the norm defined by
\[
	\|u\|_{H^{r,s}(\Sg)}^2
	:=\|u\|_{L^2(0,T;H^r(\G))}^2+\|u\|_{H^s(0,T;L^2(\G))}^2.
\]
  We clearly have
\[
	H^{0,0}(\Sg)=L^2(\Sg)=L^2(0,T;L^2(\G)).
\]
  As for the spaces with negaive exponents, we define
\[
	H^{-r,-s}(\Sg):=\(H^{r,s}_{\ ,0}(\Sg)\)',\quad r,s\ge0.
\]
  The duality paring between $H^{-r,-s}(\Sg)$ and $H^{r,s}_{\ ,0}(\Sg)$ is 
denoted by $\<\psi,u\>_{r,s}$ for $\psi\in H^{-r,-s}(\Sg)$ and 
$u\in H^{r,s}_{\ ,0}(\Sg)$.

  We define the operator $\d_{\nu_A}:H^s(\O)\to H^{s-3/2}(\G)$, $s>3/2$, as
\[
	\d_{\nu_A}u(x)
	=\frac{\d u}{\d\nu_A}(x)
	=\smij a_{ij}(x)\frac{\d u}{\d x_i}(x)\nu_j(x),
\]
  where $\nu(x)=(\nu_1(x),\dots,\nu_d(x))$ is the outward unit normal vector 
to $\G$ at $x$.
  Then for $u\in H^{r,s}(Q)$ with $r>3/2$, the trace theorem (Theorem 2.1 in 
Chapter 4 of Lions and Magenes \cite{lm}) yields
\begin{equation}\label{eq:tr}
	\frac{\d u}{\d\nu_A}\in H^{\mu,\nu}(\Sg),
	\quad\frac{\mu}{r}=\frac{\nu}{s}=\frac{r-3/2}{r}
	\quad(\mbox{$\nu=0$ if $s=0$})
\end{equation}
  and
\begin{equation}\label{eq:tr2}
	\left\|\frac{\d u}{\d\nu_A}\right\|_{H^{\mu,\nu}(\Sg)}
\le	C\|u\|_{H^{r,s}(Q)}.
\end{equation}

  In order to define the weak solution of \eqref{eq:nonhomo}, we introduce 
the {\it dual sytem};
\begin{equation}\label{eq:dual}
\begin{cases}
	\bd^{\a}v+A v=f		&\mbox{in}\ Q, \\
	v=0				&\mbox{on}\ \Sg, \\
	I^{1-\a}_{T-}v(\cdot,T)=0	&\mbox{in}\ \O,
\end{cases}
\end{equation}
  where $\bd^{\a}$ is the backward Riemann-Liouville fractional derivative, 
which is defined by
\[
	\bd^{\a}h(t)=
	-\frac{1}{\G(1-\a)}\frac{d}{dt}\int_t^T (\tau-t)^{-\a}h(\tau)d\tau,
	\quad 0<\a<1.
\]
  Moreover $I_{T-}^{\nu}$ denotes the backward integral of order $\nu$;
\[
	I^{\nu}_{T-}h(t)=
	\frac{1}{\G(\nu)}\int_t^T (\tau-t)^{\nu-1}h(\tau)d\tau,\quad \nu>0.
\]
  By the same argument as Chapter 4 in Bajlekova \cite{ba}, for any 
$f\in L^2(Q)$ there exists a unique solution $v\in H^{2,\a}(Q)$ of 
\eqref{eq:dual} such that
\begin{equation}\label{eq:vf}
	\|v\|_{H^{2,\a}(Q)}\le C\|f\|_{L^2(Q)}.
\end{equation}
  Henceforth we will denote this solution by $v_f$.
  Now we apply \eqref{eq:tr} to $v_f\in H^{2,\a}(Q)$ and obtain 
\begin{equation}\label{eq:dvf}
	\frac{\d v_f}{\d\nu_A}\in H^{1/2,\a/4}(\Sg)
	\quad\an\quad
	\left\|\frac{\d v_f}{\d\nu_A}\right\|_{H^{1/2,\a/4}(\Sg)}
	\le C\|v_f\|_{H^{2,\a}(Q)}.
\end{equation}
  Since $0<\a<1$, we have
\[
	H^{\a/4}_0(0,T;L^2(\G))=H^{\a/4}(0,T;L^2(\G))
\]
  (see \eqref{eq:Hs0}).
  That is,
\[
	\frac{\d v_f}{\d\nu_A}\in H^{1/2,\a/4}(\Sg)=H^{1/2,\a/4}_{\ ,0}(\Sg).
\]
  Now we are ready to define the weak solution of \eqref{eq:nonhomo};

\bigskip
\begin{df}
  A function $u$ is a {\it weak solution} of \eqref{eq:nonhomo} if 
\begin{equation}\label{eq:sol}
	(u,f)_{L^2(Q)}+\<g,\d_{\nu_A}v\>_{1/2,\a/4}=0
\end{equation}
  holds for any $f\in L^2(Q)$.
\end{df}

\bigskip
  The main result of this paper is as follows;

\bigskip
\begin{thm}\label{th:main}
  Let $g\in L^2(\Sg)$, then \eqref{eq:nonhomo} has a unique weak solution 
$u\in H^{1/2,\a/4}(Q)$ satisfying
\begin{equation}\label{eq:main}
	\|u\|_{H^{1/2,\a/4}(Q)}\le C\|g\|_{L^2(\Sg)}.
\end{equation}
\end{thm}

\bigskip
  Now we roughly describe the strategy of the proof.
  It is not difficult to show the unique existence of the weak solution of 
\eqref{eq:nonhomo} for $g\in L^2(\Sg)$, but the regularity of 
$H^{1/2,\a/4}(Q)$ cannot be directly deduced.
  Therefore we first show the follwoing two results;
\begin{itemize}
\item[(i)]
  Regularity of the solution for $g\in H^{-1/2,-\a/4}(\Sg)$.
\item[(ii)]
  Regularity of the solution for $g\in H^{3/2,3\a/4}(\Sg)$.
\end{itemize}
  After showing (i) and (ii), we obtain the regularity for $g\in L^2(\Sg)$ by 
interpolating the above results.

  The result for (i) can be easily shown.
  Indeed, from this definition, we can immediately deduce the following 
proposition;

\bigskip
\begin{prop}\label{prop:weak}
  Let $g\in H^{-1/2,-\a/4}(\Sg)$, then \eqref{eq:nonhomo} has a unique weak 
solution $u\in L^2(Q)$ satisfying
\begin{equation}\label{eq:weak}
	\|u\|_{L^2(Q)}\le C\|g\|_{H^{-1/2,-\a/4}(\Sg)}.
\end{equation}
\end{prop}

\bigskip
\begin{proof}
  Combining \eqref{eq:vf} and \eqref{eq:dvf}, we obtain
\begin{equation}\label{eq:vff}
	\left\|\frac{\d v_f}{\d\nu_A}\right\|_{H^{1/2,\a/4}(\Sg)}
\le	C\|f\|_{L^2(Q)}.
\end{equation}
  Thus the mapping
\[
	L^2(Q)\ni f\mapsto\frac{\d v_f}{\d\nu_A}\in H^{1/2,\a/4}(\Sg)
\]
  is bounded, and so is
\[
	L^2(Q)\ni f\mapsto-\<g,\d_{\nu_A}v_f\>\in\C.
\]
  Therefore the Riesz's representation theorem yields the unique existence of 
$u\in L^2(Q)$ such that
\[
	(u,f)_{L^2(Q)}=-\<g,\d_{\nu_A}v_f\>
\]
  holds for any $f\in L^2(Q)$.
  Thus we have proven the unique existence of weak solution.

  Moreover for any $f\in L^2(Q)$ we have
\[
	|(u,f)_{L^2(Q)}|
=	|\<g,\d_{\nu_A}v_f\>|
\le	\|g\|_{H^{-1/2,-\a/4}(\Sg)}\|\d_{\nu_A}v_f\|_{H^{1/2,\a/4}(\Sg)}
\le 	C\|g\|_{H^{-1/2,-\a/4}(\Sg)}\|f\|_{L^2(Q)},
\]
  where we have used \eqref{eq:vff} in the last inequality.
  Therefore we have \eqref{eq:weak}.
\end{proof}

\bigskip
  Thus we have proved part (i).
  In the next section, therefore, we will consider case (ii).

\section{Regular solution}\label{sec:reg}

  We first formulate the functions in $H^s(0,T)$ vanishing at $t=0$.
  Following (2.10.3.1b) in Triebel \cite{tr}, for $\R_{+}:=(0,\infty)$ we 
denote by $\wt{H}^s(\R_{+})$ the functions in $H^s(\R)$ which are identically 
zero outside $\R_{+}$;
\[
	\wt{H}^s(\R_{+}):=
	\left\{u\in H^s(\R);\ \supp u\subset\ov{\R}_{+}\right\}.
\]
  Moreover we set
\[
	\wt{H}^s_{0+}(0,T):=\{u|_{(0,T)};\ \wt{H}^s(\R_{+})\}.
\]
  Then similarly to \eqref{eq:wHs12}, we can see that such function spaces 
have a good property for interpolation;
\[
	[\wt{H}^{s_1}_{0+}(0,T),\wt{H}^{s_2}_{0+}(0,T)]_{\th}
=	\wt{H}^{(1-\th)s_1+\th s_2}_{0+}(0,T),
	\quad 0\le s_1<s_2<\infty,\ 0\le\th\le1.
\]
  We also note that we have the representation of $\wt{H}^s_{0+}(0,T)$ as
\begin{equation}\label{eq:Hs01}
	\wt{H}^s_{0+}(0,T)
	=\left\{u\in H^s(0,T);\ 
	\int_0^T|u(t)|^2\frac{dt}{t^{2s}}<\infty\right\},\quad 0\le s\le1.
\end{equation}
  In particular, for $0\le s\le1$, $s\ne1/2$, we also have
\begin{equation}\label{eq:Hs02}
	\wt{H}^s_{0+}(0,T)=
\begin{cases}
	H^s(0,T),				&0\le s<\dfrac{1}{2}, \\
	\left\{u\in H^s(0,T);\ u(0)=0\right\},	&\dfrac{1}{2}<s\le1
\end{cases}
\end{equation}
  (see \eqref{eq:Hs0} and \eqref{eq:wH0}).
  For simplicity, we set
\[
	H^{r,s}_{0+}(Q):=L^2(0,T;H^r(\O))\cap\wt{H}^s_{0+}(0,T;L^2(\O))
\]
  and
\[
	H^{r,s}_{0+}(\Sg):=L^2(0,T;H^r(\G))\cap\wt{H}^s_{0+}(0,T;L^2(\G)).
\]
  As for these spaces, we have the following properties for trace;

\bigskip
\begin{prop}\label{prop:g}
  Let $u\in H^{2,\a}_{0+}(Q)$, $0<\a<1$, then we have
\[
	u|_{\Sg}\in H^{3/2,3\a/4}_{0+}(\Sg).
\]
  Moreover the mapping 
\[
	H^{2,\a}_{0+}(Q)\ni u\mapsto u|_{\Sg}\in H^{3/2,3\a/4}_{0+}(\Sg)
\]
  is a continuous surjection.
\end{prop}

\bigskip
  For the proof of this proposition, we prepare the following lemma;

\bigskip
\begin{lem}[\bf Trace Theorem]\label{lem:trace}
  Let $X$ and $Y$ be Hilbert spaces such that $X$ is embedded to $Y$ densely 
and continuously.
  If $u\in L^2(\R_{+};X)\cap H^r(\R_{+};Y)$ with $r>1/2$, then
\[
	u^{(j)}(0)\in[Y,X]_{1-(j+1/2)/r},\quad 0\le j<r-1/2.
\]
  Moreover, the mapping 
\[
	L^2(\R_{+};X)\cap H^r(\R_{+};Y)\ni u
	\mapsto u^{(j)}(0)\in[Y,X]_{1-(j+1/2)/r}
\]
  is a continuous surjection.
\end{lem}

\bigskip
  For this lemma, see Theorem 4.2 in Chapter 1 of \cite{lm}.

\bigskip
\begin{proof}[\bf Proof of Proposition \ref{prop:g}]
  Without loss of generality, we may consider the case of
\begin{align*}
&	\O=\R_{+}^d
	=\{(x_1,\dots,x_{d-1},x_d);\ x_1,\dots,x_{d-1}\in\R,\ x_d>0\}, \\
&	\G=\R^{d-1}
	=\{(x_1,\dots,x_{d-1});\ x_1,\dots,x_{d-1}\in\R\}.
\end{align*}
  We abbreviately write $x':=(x_1,\dots,x_{d-1})$.
  Then we have
\[
	u\in H^{2,\a}_{0+}(Q)
	=L^2(0,T;H^2(\R_{+}^d))\cap\wt{H}_{0+}^{\a}(0,T;L^2(\R_{+}^d))
\]
  if and only if
\begin{align*}
	u&\in L^2\(\R_{+,x_d};H^{2,\a}_{0+}(\R^{d-1}_{x'}\tm(0,T))\)
		\cap H^2\(\R_{+,x_d};L^2(\R^{d-1}_{x'}\tm(0,T))\) \\
	&=L^2(\R_{+},H^{2,\a}_{0+}(\Sg))\cap H^2(\R_{+},L^2(\Sg)).
\end{align*}
  We apply Lemma \ref{lem:trace} as $X=H^{2,\a}_{0+}(\Sg)$ and $Y=L^2(\Sg)$.
  Then we have
\[
	u|_{\Sg}=u(0)\in[L^2(\Sg),H^{2,\a}_{0+}(\Sg)]_{3/4}
	=H^{3/2,3\a/4}_{0+}(\Sg).
\]
  The surjectivity also follows from Lemma \ref{lem:trace}.
\end{proof}

\bigskip
  By using the trace theorem stated above, problem \eqref{eq:nonhomo} with 
$g\in H^{3/2,3\a/4}_{0+}(\Sg)$ is directly reduced to the following problem 
with homogeneous boundary condition;
\begin{equation}\label{eq:homo}
\begin{cases}
	\cd^{\a}u+Au=F	&\mbox{in}\ Q, \\
	u=0		&\mbox{on}\ \Sg, \\
	u(\cdot,0)=0	&\mbox{in}\ \O.
\end{cases}
\end{equation}
  For \eqref{eq:homo}, Gorenflo, Luchko and Yamamoto \cite{gly} showed the 
$L^2$-maximal regularity.
  In their setting, the Caputo derivative $\d_t^{\a}$ equipped with the 
initial value $u(0)=0$ is formulated as an operator in $L^2(0,T)$ with its 
domain given by
\begin{equation}\label{eq:domain}
	\D(\d_t^{\a})=
\begin{cases}
	H^{\a}(0,T),					&0\le\a<1/2, \\
	\left\{u\in H^{1/2}(0,T);\ \dis\int_0^t|u(t)|^2\frac{dt}{t}<\infty
		\right\},				&\a=1/2, \\
	\left\{u\in H^{\a}(0,T);\ u(0)=0\right\},	&1/2<\a\le1.
\end{cases}
\end{equation}
  By \eqref{eq:Hs01} and \eqref{eq:Hs02}, this can be rewritten as
\begin{equation}\label{eq:domain'}
	\D(\d_t^{\a})=\wt{H}^{\a}_{0+}(0,T).
\end{equation}
  Thus we can see that if \eqref{eq:homo} has a ``solution" in 
$\D(\d_t^{\a})$, then the initial condition $u(\cdot,0)=0$ is satisfied in a 
weaker sense.
  They also revealed that the above operator $\d_t^{\a}$ is essentially 
equivalent to the Riemann-Liouville derivatives, which were already discussed 
in \cite{ba}.
  Anyway we obtain the following result;

\bigskip
\begin{lem}\label{lem:L2max}
  Let $0<\a<1$ and $F\in L^2(Q)$, then \eqref{eq:homo} has a unique solution 
$u\in H^{2,\a}_{0+}(Q)$ satisfying
\[
	\|u\|_{H^{2,\a}(Q)}\le C\|F\|_{L^2(Q)}.
\]
\end{lem}

\bigskip
  For the proof of this lemma, see Theorem 4.3 in \cite{gly}.

\bigskip
\begin{prop}\label{prop:high}
  Let $0<\a<1$ and $g\in H^{3/2,3\a/4}_{0+}(\Sg)$.
  Then problem \eqref{eq:nonhomo} has a unique solution 
$u\in H^{2,\a}_{0+}(Q)$ satisfying
\[
	\|u\|_{H^{2,\a}(Q)}\le C\|g\|_{H^{3/2,3\a/4}(\Sg)}.
\]
\end{prop}

\bigskip
\begin{proof}
  Since $g\in H^{3/2,3\a/4}_{0+}(\Sg)$, Proposition \ref{prop:g} yields that 
there exists $\tg\in H^{2,\a}_{0+}(Q)$ such that
\[
	\tg|_{\Sg}=g
	\quad\an\quad
	\|\tg\|_{H^{2,\a}(Q)}\le C\|g\|_{H^{3/2,3\a/4}(\Sg)}.
\]
  By setting $F:=-A\tg-\cd^{\a}\tg\in L^2(Q)$ and applying Lemma 
\ref{lem:L2max}, there exists a solution $w\in H^{2,\a}_{0+}(Q)$ of
\[
\begin{cases}
	\cd^{\a}w+A w=F&\mbox{in}\ Q, \\
	w=0		&\mbox{on}\ \Sg, \\
	w(\cdot,0)=0	&\mbox{in}\ \O,
\end{cases}
\]
  which satisfies
\[
	\|w\|_{H^{2,\a}(Q)}
\le 	C\|F\|_{L^2(Q)}
\le 	C\|\tg\|_{H^{2,\a}(Q)}
\le 	C\|g\|_{H^{3/2,3\a/4}(\Sg)}.
\]
  Then $u:=w+\tg$ satisfies \eqref{eq:nonhomo} and
\[
	\|u\|_{H^{2,\a}(Q)}
	\le\|w\|_{H^{2,\a}(Q)}+\|\tg\|_{H^{2,\a}(Q)}
	\le C\|g\|_{H^{3/2,3\a/4}(\Sg)}.
\]
  Thus we have completed the proof.
\end{proof}

\section{Proof of the main result}\label{sec:proof}

  In this section, we complete the proof of Theorem \ref{th:main} by 
interpolation.

\bigskip
\begin{proof}[\bf Proof of Theorem \ref{th:main}]
  Let $\pi$ be the operator which operates the boundary data $g$ to the weak
solution $u$ of \eqref{eq:nonhomo}.
  Then, by Propositions \ref{prop:weak} and \ref{prop:high}, we have
\[
	\pi\in\Lc\(H^{-1/2,-\a/4}(\Sg);L^2(Q)\)
		\cap\Lc\(H^{3/2,3\a/4}_{0+}(\Sg);H^{2,\a}(Q)\),
\]
  where $\Lc(X,Y)$ denotes the set of linear and bounded operators from $X$ to
$Y$.
  By Proposition \ref{prop:XY}, the operator $\pi$ also belongs to
\[
	\Lc\([H^{-1/2,-\a/4}(\Sg),H^{3/2,3\a/4}_{0+}(\Sg)]_{\th}
		;[L^2(Q),H^{2,\a}(Q)]_{\th}\)
\]
  for any $0\le\th\le1$.
  In particular, if we set $\th=1/4$, then
\[
	[H^{-1/2,-\a/4}(\Sg),H^{3/2,3\a/4}_{0+}(\Sg)]_{1/4}=L^2(\Sg)
	\quad\an\quad
	[L^2(Q),H^{2,\a}(Q)]_{1/4}=H^{1/2,\a/4}(Q),
\]
  and therefore we have
\[
	\pi\in\Lc(L^2(\Sg);H^{1/2,\a/4}(Q)).
\]
  Thus we have completed the proof.
\end{proof}

\appendix

\section{Interpolation}\label{sec:itp}

  Throughout this article, we often use the word ``interpolation" as a complex 
interpolation defined bellow.
  As for the detailed argument on this topic, we can refer to Triebel 
\cite{tr}, Yagi \cite{yg} and the references therein.
  On the other hand, in some classical works such as Lions and Magenes 
\cite{lm}, the ``interpolation" of two {\it Hilbert spaces} is defined as the 
domain of fractional powers of positive and self-adjoint operator.
  We will see that these two kinds of definitions coincide with each other 
(see Proposition \ref{prop:DA}).
  Therefore, we can refer to \cite{lm} and use some of their results (e.g., 
Theorem 4.2 in Chapter 1 of \cite{lm}) without any confusion.
  In this section, we recall the definition of complex interpolation of 
{\it Banach spaces} and summarize their fundamental properties.

  Let $X_i$ be a Banach space equipped with the norm $\|\cdot\|_{X_i}$ 
($i=0,1$) and suppose that $X_1$ is embedded in $X_0$ continuously and densely.
  Let $S$ be defined by
\[
	S:=\{z\in\C;\ 0<\re z<1\}.
\]
  We say that a function $F:\ov{S}\to X_0$ belongs to $\H(X_0,X_1)$ if and 
only if the following conditions (H1)-(H3) are satisfied;
\begin{itemize}
\item[\bf (H1)]
  $F$ is analytic in $S$.
\item[\bf (H2)]
  $F$ is bounded and continuous in $\ov{S}$.
\item[\bf (H3)]
  $\R\ni y\mapsto F(1+iy)\in X_1$ is bounded and continuous.
\end{itemize}
  It is known that $\H(X_0,X_1)$ is a Banach space with the norm 
$\|\cdot\|_{\H}$ given by
\[
	\|F\|_{\H}:=
	\max\(\sup_{y\in\R}\|F(iy)\|_{X_0},\sup_{y\in\R}\|F(1+iy)\|_{X_1}\),
	\quad F\in\H(X_0,X_1).
\]
  For each $0\le\th\le1$, we define the space $[X_0,X_1]_{\th}$ by
\[
	[X_0,X_1]_{\th}:=
	\{u\in X_0;\ u=F(\th)\mbox{ for some }F\in\H(X_0,X_1)\}
\]
  Moreover $[X_0,X_1]_{\th}$ is a Banach space with the norm 
$\|\cdot\|_{\th}$ defined by
\[
	\|u\|_{\th}:=
	\inf_{\begin{subarray}{c}F\in\H(X_0,X_1), \\ F(\th)=u\end{subarray}}
	\|F\|_{\H}, \quad u\in[X_0,X_1]_{\th}.
\]
  By the interpolation, we can show various kinds of 
``intermediate properties".
  For example, if a linear operator $T$ is bounded from $X_0$ into $Y_0$ and 
from $X_1$ into $Y_1$ at the same time, then we can deduce that $T$ is also a 
bounded operator from $[X_0,X_1]_{\th}$ into $[Y_0,Y_1]_{\th}$ for any 
$0<\th<1$.

\bigskip
\begin{prop}\label{prop:XY}
  Let $X_1$ (resp. $Y_1$) be embedded to $X_0$ (resp. $Y_0$) densely and 
continuously.
  Then for any $0<\th<1$,
\[
	\Lc(X_0,Y_0)\cap\Lc(X_1,Y_1)\subset\Lc([X_0,X_1]_{\th},[Y_0,Y_1]_{\th})
\]
  and we have
\[
	\|T\|_{\Lc([X_0,X_1]_{\th},[Y_0,Y_1]_{\th})}
\le	\|T\|_{\Lc(X_0,Y_0)}^{1-\th}\|T\|_{\Lc(X_1,Y_1)}^{\th},
	\quad T\in\Lc(X_0,Y_0)\cap\Lc(X_1,Y_1).
\]
\end{prop}

\bigskip
  Moreover we can also characterize the domain of fractional power of 
operators;

\bigskip
\begin{prop}\label{prop:DA}
  Let $X$ be a Hilbert space and $A:X\to X$ be a positive and self-adjoint 
operator.
  Then we have
\[
	\D(A^{\th})=[X,\D(A)]_{\th},\quad 0\le\th\le1
\]
  with isometry.
\end{prop}

\bigskip
  Here we note that $[X,\D(A)]_{\th}$ stated above coincides with 
$[\D(A),X]_{1-\th}$ in the notation by Lions and Magenes \cite{lm}.

\section{Sobolev spaces}\label{sec:sob}

  Let $\O$ be a domain of $\R^d$ with smooth boundary $\G=\d\O$.
  We dentoe by $H^s(\O)$, $s\ge0$, the space of Bessel potentials (see e.g., 
\cite{tr}).
  We can see that $H^s(\O)$, $s\ge0$, has a good property with respect to the 
interpolation;
\[
	[H^{s_1}(\O),H^{s_2}(\O)]_{\th}=H^{(1-\th)s_1+\th s_2}(\O),
	\quad 0\le s_1<s_2<\infty.
\]
  As a characterization of subspaces of $H^s(\O)$ consisting of the functions 
vanishing on $\G$, we often use $H^s_0(\O)$ (also denoted by 
$\stackrel{\circ}{H^s}(\O)$), the closure of $C_0^{\infty}(\O)$ in $H^s(\O)$.
  We note that $H^s_0(\O)$ for $0\le s\le1$ has the representation as
\begin{equation}\label{eq:Hs0}
	H^s_0(\O)=
\begin{cases}
	H^s(\O),				&0\le s\le\dfrac{1}{2}, \\
	\big\{u\in H^s(\O);\ u|_{\G}=0\big\},	&\dfrac{1}{2}<s\le1.
\end{cases}
\end{equation}
  However, we see that for these spaces, good properties concerning the 
interpolation fails for some ``singular" cases---when $s=$ integer$+1/2$.
  Indeed, according to Theorem 11.6 in Chapter 1 of \cite{lm}, the identity
\[
	[H^{s_1}_0(\O),H^{s_2}_0(\O)]_{\th}=H^{(1-\th)s_1+\th s_2}_0(\O)
\]
  is valid for the case in which $s_1,s_2,(1-\th)s_1+\th s_2\ne$integer$+1/2$.
  On the other hand, by Theorem 11.7 in Chapter 1 of \cite{lm}, if 
$(1-\th)s_1+\th s_2=\mu+1/2$ for some integer $\mu$ ($s_1,s_2$ are still 
assumed not to be integer$+1/2$), then
\begin{equation}\label{eq:0012}
	H_{00}^{\mu+1/2}(\O)
	:=[H^{s_1}_0(\O),H^{s_2}_0(\O)]_{\th}
	\varsubsetneq H^{\mu+1/2}_0(\O).
\end{equation}
  In particular, for $\mu=0$ we have
\[
	[L^2(\O),H_0^1(\O)]_{1/2}
	=H_{00}^{1/2}(\O)
	\varsubsetneq H_{0}^{1/2}(\O)
	=H^{1/2}(\O).
\]
  We note that the space $H_{00}^{\mu+1/2}(\O)$ also can be characterized 
without interpolation.
  In fact, by (11.52) in Chapter 1 of \cite{lm}, we have the repredentation as 
\begin{equation}\label{eq:m12}
	H_{00}^{\mu+1/2}(\O)
	=\left\{u\in H^{\mu+1/2}_0(\O);\ 
		\rho^{-1/2}\d_x^{\a}u\in L^2(\O)
		\ \mbox{for}\ |\a|=\mu\right\},
\end{equation}
  where $\rho:\Ov\to[0,\infty)$ is a smooth function such that
\[
	\lim_{x\to x_0}\frac{\rho(x)}{d(x,\G)}=d\ne0,\quad x_0\in\G.
\]
  Moreover by substituting $\mu=0$ in \eqref{eq:m12}, we can rewrite 
\eqref{eq:m12} as
\[
	H_{00}^{1/2}(\O)
	=\left\{u\in H^{1/2}(\O);\ 
		\int_{\O}\frac{|u(x)|^2}{\rho(x)}dx<\infty\right\}.
\]
  Thus we see that each element in $H_{00}^{1/2}(\O)$ vanishes on $\G$ in a 
certain sense.

  In the following, therefore, we introduce another formulation of functions 
in $H^s(\O)$ vanishing on $\G$, which includes $H^{\mu+1/2}_{00}(\O)$ as 
a particular case.
  We denote by $\wt{H}^s(\O)$ the subspace of $H^s(\R^d)$ consisting of the 
functions which are identically zero outside $\O$ (see (4.3.2.1b) in \cite{tr}).
  That is,
\begin{equation}\label{eq:wHs}
	\wt{H}^s(\O):=\left\{u\in H^s(\R^d);\ \supp u\subset\Ov\right\}.
\end{equation}
  By regarding elements of $\wt{H}^s(\O)$ as functions defined on $\O$, we have
\begin{equation}\label{eq:wH0}
	\wt{H}^s(\O)=H^s_0(\O),\quad s\ge0,\ s\ne\mbox{integer}+\frac{1}{2}.
\end{equation}
  Moreover, we have the following good property for interpolation (see 
Corollary 1.6 in Chapter 3 of Strichartz \cite{str});
\begin{equation}\label{eq:wHs12}
	[\wt{H}^{s_1}(\O),\wt{H}^{s_2}(\O)]_{\th}
	=\wt{H}^{(1-\th)s_1+\th s_2}(\O),
	\quad 0\le s_1<s_2<\infty,\ 0\le\th\le1.
\end{equation}
  By comparing \eqref{eq:0012} with \eqref{eq:wHs12} and noting 
\eqref{eq:wH0}, we have
\[
	H^{\mu+1/2}_{00}(\O)=\wt{H}^{\mu+1/2}(\O),\quad \mu\in\N_0.
\]
  This identity also can be verified by comparing the representations 
\eqref{eq:m12} and (2.4.2.7) in \cite{tr} with $p=2$.
  Thus the subspace $\wt{H}^s(\O)$ introduced in \eqref{eq:wHs} is more 
appropriate than $H_0^s(\O)$ when we deal with interpolation.


\end{document}